\documentclass[12pt, twoside, eqno]{article}
\usepackage{latexsym}
\usepackage{amssymb}
\usepackage{amsfonts}
\usepackage{graphicx}
\begin{document}
\noindent {\bf \Large On involution
$le$-semigroups}\bigskip

\medskip

\noindent{\bf Niovi Kehayopulu}\bigskip

\bigskip{\small

\noindent{\bf Abstract.} We deal with involution ordered semigroups 
possessing a greatest element, we introduce the concepts of 
$*$-regularity, $*$-intra-regularity, $*$-bi-ideal element and 
$*$-quasi-ideal element in this type of semigroups and, using the 
right and left ideal elements, we give relations between the 
regularity and $*$-regularity, between intra-regularity and 
$*$-intra-regularity. Finally, we prove that in an involution 
$*$-regular $\vee e$-semigroup every $*$-bi-ideal element can be 
considered as a product of a right and a left ideal element, we 
describe the form of the filter generated by an element of an 
involution $*$-intra-regular $poe$-semigroup $S$, showing that every 
$\cal N$-class of $S$ has a greatest element.\medskip

\noindent{\bf 2010 AMS Subject Classification:} 06F05.\smallskip

\noindent{\bf Keywords:} $poe$ ($le$)-semigroup; involution; right 
(left) ideal element, bi-ideal element; $*$-right ($*$-left) ideal 
element; $*$-bi-ideal element; regular; intra-regular; $*$-regular; 
$*$-intra-regular. }
\section{Introduction} Semigroups with involution (called as 
$*$-semigroups as well) have been widely studied.
Involution regular semigroups have been studied, among others, in 
[2--3, 5]. Wu, Chong-Yih recently introduced and studied the 
involution ordered semigroups ($po$-semigroups) in Thai J. Math. But 
we should mention that if we have an ordered semigroup and consider 
the (involution) ``$*$" as the identity mapping, we do not have an 
involution ordered semigroup in general. So we cannot say that the 
results on involution ordered semigroups generalize results on 
ordered semigroups.
In the present paper we deal with involution ordered semigroups 
possessing a greatest element ($poe$, $le$-semigroups) using elements 
instead of sets. We introduce the concepts of involution $*$-regular 
and involution $*$-intra-regular $poe$-semigroups and we prove, among 
others, the following: If an involution $poe$-semigroup $S$ is 
$*$-regular (resp. $*$-intra-regular), then for every left ideal 
element $a$ and every right ideal element $b$ of $S$ such that 
$a\wedge b$ exists, we have $a\wedge b\le a^* b^*$ (resp. $a\wedge 
b\le b^* a^*$). Conversely, if $S$ is an involution $le$-semigroup 
such that $a\wedge b\le b^* a^*$ (resp. $a\wedge b\le a^* b^*$) for 
every left ideal element $a$ and every right ideal element $b$ of 
$S$, then $S$ is regular (resp. intra-regular). An involution 
$le$-semigroup $S$ is $*$-regular if and only if (1) the right and 
the left ideal elements of $S$ are idempotent and for any right ideal 
element $a$ and any left ideal element $b$ of $S$, the product $ab$ 
is a quasi-ideal element of $S$ and (2) if $a\le r(a^2)$ and $a\le 
l(a^2)$ for any $a\in S$. An involution $poe$-semigroup $S$ is 
$*$-intra-regular if and only if the ideal elements of $S$ are 
$*$-semiprime and this is equivalent to saying that $N(x)=\{y\in S 
\mid x\le ey^*e\}$ for every $x\in S$ and so each class $(x)_{\cal 
N}$ has a greatest element.
\section{Prerequisites} An ordered groupoid (: $po$-groupoid) is a 
groupoid $S$ with an order relation ``$
\le$" on $S$ such that $a\le b$ implies $ac\le bc$ and $ca\le cb$ for 
every $c\in S$. A $\vee$-groupoid is a $po$-groupoid $S$ at the same 
time a semilattice under $\vee$ such that $(a\vee b)c=ab\vee ac$ and 
$a(b\vee c)=ab\vee ac$ for all $a,b,c\in S$. A $\vee$-groupoid which 
is at the same time a lattice is called an $l$-groupoid [1]. If the 
multiplication on a $po$-groupoid $S$ is associative, then $S$ is 
called a $po$-semigroup. By a $poe$-groupoid we mean a $po$-groupoid 
possessing a greatest element $e$ (that is, $e\ge a$ $\forall$ $a\in 
S$). An element $a$ of a  $po$-groupoid $S$ is called idempotent if 
$a^2=a$; it is called a left (resp. right) ideal element if $xa\le a$ 
(resp. $ax\le a$) for every $x\in S$ (cf. also [1]). In a 
$poe$-groupoid,
an element $a$ is a left (resp. right) ideal element if and only if 
$ea\le a$ (resp. $ae\le a$). The element $ea$ (resp. $ae$) is a left 
(resp. right) ideal element of $S$; and the element $eae$ is an ideal 
element of $S$. We denote by $l(a)$ and $r(a)$ the left and the right 
ideal element of $S$, respectively, generated by $a$ $(a\in S$), and 
one can easily prove that $l(a)=a\vee ea$ and $r(a)=a\vee ae$. 
Moreover, for every $a,b\in S$, we have $a\le l(a)$ and $a\le r(a)$;
$a\le b$ implies $l(a)\le l(b)$ and $r(a)\le r(b)$; $l((l(a))=l(a)$, 
$r((r(a))=r(a)$; $l((r(a))=r((l(a))$. An element $a$ of a 
$po$-groupoid $S$ is called {\it semiprime} if for any $t\in S$ such 
that $t^2\le a$, we have $t\le a$. An element $a$ of a $poe$-groupoid 
$S$ is called a quasi-ideal element of $S$ if the element $ae\wedge 
ea$ exists in $S$ and we have $ae\wedge ea\le a$; in particular if 
the $S$ is a $poe$-semigroup, then the element $a$ is called a 
bi-ideal element of $S$ if $aea\le a$. A subsemigroup $F$ of $S$ is 
called a {\it filter} of $S$ if (1) $a,b\in S$, $ab\in F$ implies 
$a\in F$ and $b\in F$ and
(2) if $a\in F$ and $S\ni b\ge a$, then $b\in F$. For an element $x$ 
of $S$, we denote by $N(x)$ the filter of $S$ generated by $x$.
For a $po$-groupoid $S$ and a subset $H$ of $S$ we denote by $(H]$ 
the subset of $S$ defined by $(H]:=\{t\in S \mid t\le a \mbox { for 
some } a\in H\}.$ An element $a$ of a $po$-semigroup $S$ is called 
regular if there exists an element $x\in S$ such that $a\le axa$, 
that is if $a\in (aSa]$. An element $a$ of a $po$-semigroup $S$ is 
called intra-regular if there exist $x,y\in S$ such that $a\le xa^2 
y$, that is if $a\in (Sa^2 S]$. A $po$-semigroup $S$ is regular 
(resp. intra-regular) if every element of $S$ is so. Recall that $S$ 
is regular (resp. intra-regular) if and only if $A\subseteq (ASA]$ 
(resp. $A\subseteq (AS^2 A])$ for any subset $A$ of $S$. A 
$poe$-semigroup $S$ is regular (resp. intra-regular) if and only if 
$a\le aea$ (resp. $a\le ea^2 e)$ for every $a\in S$.

\section{Main results}{\bf Definition 1.} [6] An ordered groupoid 
$(S,\cdot,\le)$ is called {\it involution ordered groupoid} (: {\it 
involution po-groupoid}) if there exists a unary operation ``$*$" on 
$S$ such that

(1) $(a^*)^*=a$ and $(ab)^*=b^* a^*$ for every $a,b\in S$ and

(2) if $a\le b$, then $a^*\le b^*$.\\
In an involution $poe$-groupoid, we clearly have $e^*=e$.\smallskip

\noindent{\bf Example 2.} (cf. also [4; the Example]) The set 
$S=\{a,b,c,d,e\}$ with the multiplication ``$\cdot$", the involution 
``$*$", and the figure below is an example of an involution 
$po$-semigroup. \begin{center}
$\begin{array}{*{20}{c}}
{\cdot}&\vline& a&\vline& b&\vline& c&\vline& d&\vline& f\\
\hline
a&\vline& a&\vline& a&\vline& a&\vline& a&\vline& a\\
\hline
b&\vline& a&\vline& b&\vline& a&\vline& d&\vline& a\\
\hline
c&\vline& a&\vline& f&\vline& c&\vline& c&\vline& f\\
\hline
d&\vline& a&\vline& b&\vline& d&\vline& d&\vline& b\\
\hline
f&\vline& a&\vline& f&\vline& a&\vline& c&\vline& a
\end{array}$
\end{center}
$$a^*=a,\, b^*=c,\, d^*=d,\, c^*=b,\, f^*=f$$
\begin{center}
\begin{figure}[h]
\centering
\includegraphics[width=0.31\textwidth]{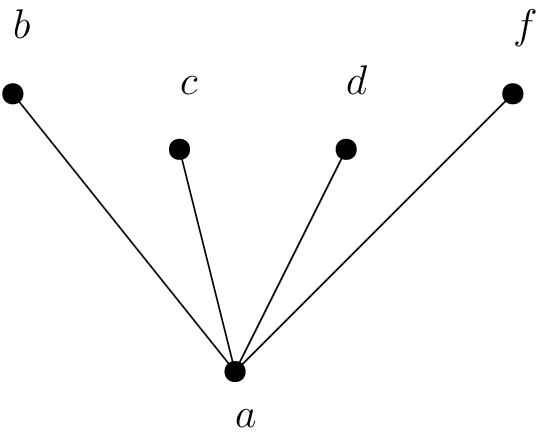}
\end{figure}\end{center}
\noindent{\bf Definition 3.} Let $S$ be an involution $poe$-groupoid. 
An element $a$ of $S$ is called a {\it $*$-right} (resp. {\it 
$*$-left}) {\it ideal element} of $S$ if $a^*e\le a$ (resp. $ea^*\le 
a)$; if $S$ is in addition semilattice under $\wedge$, then $a$ is 
called a {\it $*$-quasi-ideal element} of $S$ if $a^*e\wedge ea^*\le 
a$.
If $S$ is an involution $poe$-semigroup, an element $a$ of $S$ is 
called a {\it $*$-bi-ideal element} of $S$ if $a^*ea^*\le 
a$.\smallskip

\noindent{\bf Proposition 4.} {\it If S is an involution poe-groupoid 
at the same time semilattice under $\wedge$, then every $*$-right} 
({\it resp. $*$-left}) {\it ideal element a of S is a $*$-quasi-ideal 
element of S.
If S is an involution poe-semigroup which is also semilattice under 
$\wedge$, then every $*$-quasi-ideal element of S is a $*$-bi-ideal 
element of S}.\smallskip

\noindent{\bf Proof.} Let $a$ be a $*$-right ideal element of $S$, 
that is $a^*e\le a$. Then $a^*e\wedge ea^*\le a^*e\le a$, so $a$ is a 
$*$-quasi-ideal element of $S$. If $a$ is a $*$-quasi-ideal element 
of $S$, then $a^*ea^*\le a^*e\wedge ea^*\le a$, so $a^*$ is a 
$*$-bi-ideal element of $S$.\smallskip

\noindent{\bf Proposition 5.} {\it In an involution po-groupoid which 
is at the same time semilattice under $\vee$} ({\it resp. under} 
$\wedge$), {\it we have}$$(a\vee b)^*=a^*\vee b^* \mbox { ({\it 
resp.} } (a\wedge b)^*=a^*\wedge b^*).$${\bf Proof.} Since $a\vee 
b\ge a, b$, we have $(a\vee b)^*\ge a^*, b^*$. Let now $t\in S$ such 
that $t\ge a^*$ and $t\ge b^*$. Then $(a\vee b)^*\le t$. Indeed: 
Since $t^*\ge (a^*)^*=a$ and $t^*\ge (b^*)^*=b$, we have $t^*\ge 
a\vee b$, then $t=(t^*)^*\ge (a\vee b)^*$ and so $(a\vee b)^*=a^*\vee 
b^*$. The dual case can be proved in a similar way. 
$\hfill\Box$\smallskip

\noindent{\bf Proposition 6.} {\it In an involution $\vee 
e$-semigroup S, for every $a\in S$, we have

$al(a^*)=r(a)a^*$ and $r(a^*)a=a^*l(a)$.} \smallskip

\noindent{\bf Proof.} Indeed, we have

$al(a^*)=a(ea^*\vee a^*)=aea^*\vee aa^*=(ae\vee a)a^*=r(a)a^*$ and

$r(a^*)a=(a^* e\vee a^*)a=a^* ea\vee a^* a=a^*(ea\vee a)=a^*l(a)$.
 $\hfill\Box$

\noindent{\bf Proposition 7.} {\it Let S be an involution 
poe-groupoid. Then a is a left} ({\it resp. right}) {\it ideal 
element of S if and only if $a^*$ is a right} ({\it resp. left}) {\it 
ideal element of S; if S is at the same time semilattice under 
$\wedge$, then a is a quasi-ideal element of S if and only if $a^*$ 
is so.
If S is an involution poe-semigroup, then a is a bi-ideal element of 
S if and only if $a^*$ is so.}\smallskip

\noindent{\bf Proof.} Let $a$ be a left ideal element of $S$. Since 
$ea\le a$, we have $(ea)^*\le a^*$, then $a^* e=a^* e^*=(ea)^*\le 
a^*$, so $a^*$ is a right ideal element of $S$. If $a$ is a right 
ideal element of $S$, then $ae\le a$, so $ea^*\le a^*$, and $a^*$ is 
a left ideal element of $S$. From this follows that if $a^*$ is a 
right (resp. left) ideal element of $S$, then the element 
$(a^*)^*(=a)$ is a left (resp. right) ideal element of $S$, so $a$ is 
a left (right) ideal element of $S$ if and only if $a^*$ is a right 
(left) ideal element of $S$. Let now $S$ be a semilattice under 
$\wedge$ and $a$ a quasi-ideal element of $S$. Since $ae\wedge ea\le 
a$, we have $(ae\wedge ea)^*\le a^*$. Then, by Proposition 5, 
$(ae)^*\wedge (ea)^*\le a^*$, then $ea^*\wedge a^*e\le a^*$ and so 
$a^*$ is a quasi-ideal element of $S$. Conversely, if $a^*$ is a 
quasi-ideal element of $S$, then $a^*e\wedge ea^*\le a^*$, then 
$(a^*e\wedge ea^*)^*\le a$, $(a^*e)^*\wedge (ea^*)^*\le a$, and 
$ea\wedge ae\le a$, so $a$ is a quasi-ideal element of $S$. If $S$ is 
an involution $poe$-semigroup then we have $aea\le a$ if and only if 
$(aea)^*\le a^*$ and this is equivalent to $a^*ea^*\le a^*$, so $a^*$ 
is a bi-ideal element of $S$ if and only if $a^*$ is so. 
$\hfill\Box$\smallskip

\noindent{\bf Proposition 8.} {\it If S is an involution poe-groupoid 
then for every left ideal element a and every right ideal element b 
of S such that $(a^*\wedge b^*)e\wedge e(a^*\wedge b^*)$ exists, 
$a^*\wedge b^*$ is a quasi-ideal element of S.
If S is involution regular poe-semigroup and a a left} ({\it resp. 
right}) {\it ideal element of S, then the element $a^*$ is 
idempotent}.\smallskip

\noindent{\bf Proof.} Let $a$ be a left ideal element and $b$ a right 
ideal element of $S$ such that $(a^*\wedge b^*)e\wedge e(a^*\wedge 
b^*)$ exists. Since $a^*$ is a right and $b^*$ a left ideal element 
of $S$, we have$$(a^*\wedge b^*)e\wedge e(a^*\wedge b^*)\le 
a^*e\wedge eb^*\le a^*\wedge b^*,$$so $a^*\wedge b^*$ is a 
quasi-ideal element of $S$.
Let now $S$ be an involution regular $poe$-semigroup. If $a$ a left 
ideal element of $S$, then $a^*\le (a^*e)a^*\le a^* a^*\le a^*e\le 
a^*$, so $a^* a^*=a^*$. If $a$ a right ideal element of $S$, then we 
have $a^*\le a^*(ea^*)\le a^* a^*\le ea^*\le a^*$, then $a^* 
a^*=a^*$. $\hfill\Box$\smallskip

\noindent{\bf Proposition 9.} {\it Let S be an involution 
poe-semigroup. If a is a right ideal element of S and $b\in S$, then 
the element $(ab)^*$ is a bi-ideal element of S. If b is a right 
ideal element and $a\in S$, then $(ab)^*$ is a bi-ideal element of 
S}.\smallskip

\noindent{\bf Proof.} Let $a$ be a right ideal element and $b\in S$. 
Then have$$(ab)^*e(ab)^*= b^*a^*eb^*a^*\le b^*(ea^*)\le 
b^*a^*=(ab)^*.$$If $b$ is a right ideal element of $S$ and $a\in S$, 
then
$$(ab)^*e(ab)^*= b^*a^*eb^*a^*\le (eb^*)a^*\le b^*a^*=(ab)^*$$and the 
proof is completed.

Another proof is to show that if $a$ is a right ideal element of $S$ 
and $b\in S$ (or if $b$ is a right ideal element of $S$ and $a\in 
S$), then $ab$ is a bi-ideal element of $S$. Then, by Proposition 7, 
$(ab)^*$ is a bi-ideal element of $S$. $\hfill\Box$ \smallskip

\noindent{\bf Definition 10.} Let $S$ be an involution $po$-groupoid. 
An element $a$ of $S$ is called {\it $*$-semiprime} if for any $t\in 
S$ such that $t^*t^*\le a$, we have $t\le a$.\smallskip

\noindent{\bf Proposition 11.} {\it Let S be an involution 
poe-semigroup. Then we have the following:
\begin{enumerate}
\item[$(1)$] if $x\in I(x^*x^*)$ for every $x\in S$, then the 
ideal elements of S are $*$-semiprime;
\item[$(2)$] if $x\in I(x^2)$ for every $x\in S$, then the ideal 
elements of S are semiprime.\end{enumerate}}
\noindent{\bf Proof.} (1) Let $a$ be an ideal of $S$ and $t\in S$ 
such that $t^*t^*\le a$. Then $t\le a$. Indeed: Since $t\in 
I(t^*t^*)$, we have $t\le t^*t^*$ or $t\le et^*t^*$ or $t\le t^*t^*e$ 
or $t\le et^*t^*e$. If $t\le t^*t^*$, then $t\le a$; if $t\le 
et^*t^*$, then $t\le ea\le a$; if $t\le t^*t^*e$, then $t\le ae\le 
a$; if $t\le et^*t^*e$, then $t\le eae\le a$.\\(2) Let $a$ be an 
ideal element of $S$ and $t\in S$ such that $t^2\le a$. Then $t\le 
a$. Indeed: Since $t\in I(t^2)$, we have $t\le t^2$ or $t\le et^2$ or 
$t\le t^2e$ or $t\le et^2e$. In each case $t\le a$. 
$\hfill\Box$\smallskip

\noindent{\bf Definition 12.} An involution $poe$-semigroup $S$ is 
called {\it $*$-regular} if, for every $a\in S$, we have $a\le a^* 
ea^*.$\smallskip

\noindent{\bf Theorem 13.} {\it Let S be an involution poe-semigroup. 
If S is $*$-regular, then for every left ideal element a and every 
element b of S} ({\it or any right ideal element b and any $a\in S$}) 
{\it such that $a\wedge b$ exists, we have} $a\wedge b\le a^* b^*$. 
{\it For the converse statement, suppose S is an} ({\it involution}) 
{\it le-semigroup and for every left ideal element a and every right 
ideal element b of S, we have $a\wedge b\le b^*a^*,$ then S is 
regular.}\smallskip

\noindent{\bf Proof.} $\Longrightarrow$. Let $a$ be a left ideal 
element of $S$ and $b\in S$ such that $a\wedge b$ exists. Since $S$ 
is $*$-regular and $a^*$ a right ideal element of $S$, we have

$a\wedge b\le (a\wedge b)^*e(a\wedge b)^*=(a^*\wedge b^*)e(a^*\wedge 
b^*)\le (a^*e)b^*\le a^*b^*.$\\If $b$ is a right ideal element of $S$ 
and $a\in S$, then $a\wedge b\le a^*(eb^*)\le 
a^*b^*.$\\$\Longleftarrow$. Let $a\in S$. Since $r(a)$ (resp. $l(a)$) 
is a right (resp. left) ideal element of $S$, by hypothesis, we 
have\begin{eqnarray*}a&\le&r(a)\wedge l(a)\le (l(a))^*(r(a))^*=(a\vee 
ea)^*(a\vee ae)^*\\&=&(a^*\vee a^*e)(a^*\vee ea^*)=a^*a^*\vee (a^* 
e)a^*\vee a^*(ea^*)\vee (a^* e)(ea^*)\\&=&a^*a^*\vee 
a^*ea^*.\end{eqnarray*}Then we have \begin{eqnarray*}a^*&\le& 
(a^*a^*\vee a^*ea^*)^*=(a^*a^*)^*\vee (a^*ea^*)^*\mbox { (by 
Proposition 5)}\\&=&(a^*)^*(a^*)^*\vee (a^*)^*e^*(a^*)^*=a^2\vee 
aea,\end{eqnarray*}from which$$a^*a^*\le (a^2\vee aea)(a^2\vee 
aea)=a^4\vee aea^3\vee a^3ea\vee aea^2ea\le aea,$$and
\begin{eqnarray*}a^*ea^*&\le&(a^2\vee aea)e(a^2\vee 
aea)\\&=&a^2ea^2\vee aeaea^2\vee a^2eaea\vee aeaeaea\le 
aea.\end{eqnarray*}Thus we get $a\le aea$, and $S$ is 
regular.$\hfill\Box$\smallskip

\noindent{\bf Proposition 14.} {\it If S is an involution $*$-regular 
$\vee e$-semigroup then, for every $a\in S$, we have $a\le r(a^*)$ 
and $a\le l(a^*)$}.\smallskip

\noindent{\bf Proof.} Let $a\in S$. Since $S$ is $*$-regular, we 
have$$a\le a^*(ea^*)\le a^*e\le a^*e\vee a^*=r(a^*),$$and $a\le 
(a^*e)a^*\le ea^*\le l(a^*)$. $\hfill\Box$\smallskip

\noindent{\bf Proposition 15.} {\it If S is an involution $*$-regular 
poe-semigroup and a a left} ({\it or right}) {\it ideal element of S, 
then $a=a^*$; if b a bi-ideal element of S, then $b=b^*$.}\smallskip

\noindent{\bf Proof.} Let $a$ be a left ideal element of $S$. Since 
$S$ is $*$-regular and $a^*$ a right ideal element of $S$, we have 
$a\le (a^*e)a^*\le a^*a^*\le a^*e\le a^*$, then $a^*\le (a^*)^*=a$, 
thus we have $a=a^*$. If $a$ is a right ideal element of $S$, then 
$a\le
a^*(ea^*)\le a^*a^*\le ea^*\le a^*$, then $a^*\le a$, so again 
$a=a^*$. If $b$ is a bi-ideal element of $S$, then we have $beb\le 
b$, then $b^*eb^*\le b^*$; since $S$ is $*$-regular, we have $b\le 
b^*eb^*$, thus $b\le b^*$, then $b^*\le b$, and so $b^*=b$. 
$\hfill\Box$\smallskip

\noindent{\bf Proposition 16.} {\it Let S be an involution 
$*$-regular poe-semigroup which is also semilattice under $\wedge$. 
If a is a right ideal element and b a left ideal element of S, then 
the product $ab$ is a quasi-ideal element of S}.\smallskip

\noindent{\bf Proof.} Let $a$ be a right ideal element and $b$ a left 
ideal element of $S$. Since $S$ is $*$-regular, by Proposition 15, we 
have $a^*=a$ and $b^*=b$, we also have\begin{eqnarray*}a\wedge b&\le& 
(a\wedge b)^*e(a\wedge b)^*=(a^*\wedge b^*)e(a^*\wedge b^*)\le 
(a^*e)b^*\\&=&(ae)b\le ab\le ae\wedge eb\le a\wedge 
b,\end{eqnarray*}so $a\wedge b=ab$. On the other hand, $(a\wedge 
b)e\wedge e(a\wedge b)\le ae\wedge eb\le a\wedge b$, that is, 
$a\wedge b$ is a quasi-ideal element of $S$; then $ab$ is a 
quasi-ideal element element of $S$ as well.$\hfill\Box$\smallskip

\noindent{\bf Proposition 17.} {\it An involution $*$-regular 
poe-semigroup S is regular. In regular poe-semigroups, the right and 
the left ideal elements are idempotent.}\smallskip

\noindent{\bf Proof.} Let $a\in S$. Since $S$ is $*$-regular, we have 
$a\le a^*ea^*$. Then we have $a^*\le aea$, then $a\le (aea)e(aea)\le 
aea$ and so $S$ is regular. Let now $a$ be a right ideal element of 
$S$. Since $S$ is regular, we have $a\le (ae)a\le a^2\le ae\le a$, so 
$a$ is idempotent. If $a$ is a left ideal element of $S$, then $a\le 
a(ea)\le a^2\le ea\le a$, so $a^2=a$.\smallskip

\noindent{\bf Proposition 18.} {\it If S is an involution 
le-semigroup such that\begin{enumerate}
\item[$(1)$] $a\le r(a^*)$ and $a\le l(a^*)$ for every $a\in S$ and
\item[$(2)$] For every right ideal element a and every left ideal 
element b of S, we have

    $a^2=a$, $b^2=b$ and ab is a quasi-ideal element of 
S,\end{enumerate}then S is $*$-regular.}

\noindent{\bf Proof.} Let $a\in S$. Since $r(a^*)$ is a right ideal 
element of $S$, we have\begin{eqnarray*}a&\le& 
r(a^*)=(r(a^*))^2=(a^*e\vee a^*)(a^*e\vee a^*)\\
&=&a^*ea^*e\vee a^*a^*e\vee a^*ea^*\vee a^*a^*\\&\le& 
a^*e,\end{eqnarray*}similarly $a\le ea^*$, hence $a\le a^*e\wedge 
ea^*$. Since $a^*e$ is a right ideal element and $ea^*$ a left ideal 
element element of $S$, we have $$a^*e\wedge ea^*=(a^*e)^2\wedge 
(ea^*)^2=(a^*ea^*)e\wedge e(a^*ea^*).$$ Since $e=e^*$, we have 
$a^*ea^*=(a^*e)(ea^*).$ Since $a^*e$ is a right ideal element and 
$ea^*$ a left ideal element of $S$, by hypothesis, $(a^*e)(ea^*)$ is 
a quasi-ideal element of $S$. Then $a^*ea^*$ is a quasi-ideal element 
of $S$, that is, $(a^*ea^*)e\wedge e(a^*ea^*)\le a^*ea^*$. Then we 
have $a\le a^*ea^*$, and $S$ is $*$-regular. $\hfill\Box$.\\
By Propositions 14, 16--18, we have the following theorem\smallskip

\noindent{\bf Theorem 19.} {\it An involution le-semigroup is 
$*$-regular if and only if conditions $(1)$ and $(2)$ of Proposition 
$18$ are satisfied}.\smallskip

\noindent{\bf Theorem 20.} {\it Let S is a $*$-regular involution 
$\vee e$-semigroup and b a $*$-bi-ideal element of S. Then there 
exist a right ideal element x and a left ideal element y of S such 
that $b=xy$}.\smallskip

\noindent{\bf Proof.} Since $S$ is $*$-regular and $b$ a $*$-bi-ideal 
element of $S$, we have $b=b^*eb^*$. Then we have
\begin{eqnarray*}r(b^*)l(b^*)&=&(b^*\vee b^*e)(b^*\vee 
eb^*)=b^*b^*\vee b^*eb^*\vee b^*e^2b^*\\&=&b^*b^*\vee 
b^*eb^*=b^*b^*\vee b.\end{eqnarray*}
Moreover, $b^*b^*=(beb)(beb)\le beb=b^*eb^*=b$ (by Prop. 15). So we 
have $r(b^*)l(b^*)=b$, where $r(b^*)$ is a right ideal element and 
$l(b^*)$ a left ideal element of $S$.\smallskip

\noindent{\bf Definition 21.} An involution $poe$-semigroup $S$ is 
called {\it $*$-intra-regular} if, for every $a\in S$, we have $a\le 
ea^*a^*e.$ \smallskip

\noindent{\bf Theorem 22.} {\it Let S be an involution poe-semigroup. 
If S is $*$-intra-regular, then for every left ideal element a and 
every right ideal element b of S such that $a\wedge b$ exists, we 
have $a\wedge b\le b^* a^*.$ ``Conversely", if S be an involution 
$le$-semigroup such that for every left ideal element a and every 
right ideal element b of S, we have $a\wedge b\le a^*b^*,$ then S is 
intra-regular.}\smallskip

\noindent{\bf Proof.} $\Longrightarrow$. Let $a$ be a left ideal 
element and $b$ a right ideal element of $S$ such that $a\wedge b$ 
exists. Since $S$ is $*$-intra-regular, we have$$a\wedge b\le 
e(a\wedge b)^*(a\wedge b)^*e\le (eb^*)(a^*e)\le 
b^*a^*.$$$\Longleftarrow$. Let $a\in S$. By hypothesis, we 
have\begin{eqnarray*}a&\le& r(a)\wedge l(a)\le 
(r(a))^*(l(a))^*=(a\vee ae)^*(a\vee ea)^*\\&=&(a^*\vee ea^*)(a^*\vee 
a^*e) =a^*a^*\vee ea^*a^*\vee a^*a^*e\vee 
ea^*a^*e.\end{eqnarray*}Then we have
\begin{eqnarray*} a^*&\le&(a^*a^*\vee ea^*a^*\vee a^*a^*e\vee 
ea^*a^*e)^*\\&=&(a^*a^*)^*\vee (ea^*a^*)^*\vee (a^*a^*e)^*\vee 
(ea^*a^*e)^* \mbox { (by Proposition 5)}\\&=&a^2\vee a^2e\vee 
ea^2\vee ea^2e.\end{eqnarray*}Then $$a^*a^*\le (a^2\vee a^2e\vee 
ea^2\vee ea^2e)(a^2\vee a^2e\vee ea^2\vee ea^2e)\le ea^2e\vee 
a^2e,$$
$$ea^*a^*\le e(ea^2e\vee a^2e)=e^2a^2e\vee ea^2e=ea^2e,$$
$$a^*a^*e\le (ea^2e\vee a^2e)e=ea^2e^2\vee a^2e^2\le ea^2e\vee 
a^2e,$$$$ea^*a^*e=(ea^*a^*)e\le (ea^2e)e=ea^2e^2\le ea^2e.$$Hence we 
get $a\le ea^2e\vee a^2e$, then$$a^2 e\le a(ea^2e\vee 
a^2e)e=aea^2e\vee a^3e^2\le ea^2e,$$and $a\le ea^2e$, so $S$ is 
intra-regular.$\hfill\Box$\smallskip

\noindent{\bf Proposition 23.} {\it Let S be an involution 
$poe$-semigroup. If S in $*$-intra-regular, then $eabe=eb^* a^* e$ 
for every $a,b\in S$}.\smallskip

\noindent{\bf Proof.} Let $a,b\in S$. Since $S$ is $*$-intra-regular, 
we have$$ab\le e(ab)^*(ab)^*e=eb^*a^*b^*a^* e\le eb^*a^* e,$$then 
$eabe\le eb^*a^* e$. Again since $S$ is $*$-intra-regular, we 
have$$b^*a^*\le e(b^*a^*)^*(b^*a^*)^*e=e(ab)(ab)e\le eabe,$$then 
$eb^*a^* e\le eabe$. Thus we have $eabe=eb^* a^* e$. $\hfill\Box$
\smallskip

\noindent{\bf Proposition 24.} {\it If an involution poe-semigroup S 
is $*$-intra-regular, then the ideal elements of S are $*$-semiprime. 
``Conversely", if the ideal elements of S are $*$-semiprime, then $S$ 
is intra-regular}.\smallskip

\noindent{\bf Proof.} $\Longrightarrow$. Let $a$ be an ideal element 
of $S$ and $t\in S$ such that $t^*t^*\le a$. Since $S$ is 
$*$-intra-regular, we have $t\le et^*t^*e\le eae\le a$.\\
$\Longleftarrow$. Let $a\in S$. Since the element $ea^*a^*e$ is an 
ideal element of $S$, by hypothesis, it is semiprime; and since 
$(aa)^*(aa)^*=a^*a^*a^*a^*\le ea^*a^*e$, we have $aa\le ea^*a^*e$, 
then $a^*a^*=(aa)^*\le eaae$. Since $eaae$ is an ideal element of 
$S$, it is $*$-semiprime. Thus we have $a\le eaae=ea^2e$, and $S$ is 
intra-regular. $\hfill\Box$\smallskip

\noindent{\bf Proposition 25.} {\it If an involution poe-semigroup S 
is $*$-regular} ({\it resp. $*$-intra-regular}), {\it then it is 
regular} ({\it resp. intra-regular}).\smallskip

\noindent{\bf Proof.} Let $a\in S$. If $S$ is $*$-regular, we have 
$a\le a^*ea^*$, then $a^*\le aea$, so
$a\le (aea)e(aea)\le aea$, and $S$ is regular. Similarly, if $S$ is
$*$-intra-regular, then it is intra-regular. $\hfill\Box$\smallskip

\noindent{\bf Theorem 26.} {\it An involution poe-semigroup S is 
$*$-intra-regular if and only if, for every $x\in S$, we 
have$$N(x)=\{y\in S \mid x\le ey^*e\}.$$}
{\bf Proof.} $\Longrightarrow$. Let $x\in S$ and $T:=\{y\in S \mid 
x\le ey^*e\}.$ $T$ is a nonempty subset of $S$. Indeed, since $S$ is 
$*$-intra-regular, we have $x\le ex^*x^*e$, thus $x^*x^*\in T$.\\Let 
$a,b\in T$. Since $x\le ea^*e$ and $x\le eb^*e$, by hypothesis, we 
have\begin{eqnarray*}x&\le& ex^*x^*e\le 
e(eb^*e)^*(ea^*e)^*e=e(ebe)(eae)e\\&\le&e(bea)e\le ee(bea)^*(bea)^*e 
\mbox { (since } S \mbox { is 
$*$-intra-regular)}\\&\le&e(a^*eb^*)(a^*eb^*)e\le 
e(b^*a^*)e\\&=&e(ab)^*e,\end{eqnarray*}so $ab\in T$. Let $a,b\in S$ 
such that $ab\in T$. Then $x\le e(ab)^*e=eb^*a^*e\le ea^*e, eb^*e,$ 
then $a,b\in T$. Let $a\in T$ and $S\ni b\ge a$. Then $x\le ea^*e\le 
eb^*e$, so $b\in T$.\\Let now $F$ be a filter of $S$ such that $x\in 
F$ and $a\in T$. Since $S$ is $*$-intra-regular, we have $F\ni x\le 
ea^*e\le ee(a^*)^*(a^*)^*ee=eaae\le eae,$ then $eae\in F$, and $a\in 
F$.\\$\Longleftarrow$. Let $x\in S$. Since $N(x)$ is a subsemigroup 
of $S$ containing $x$, we have $xx\in N(x)$, then $x\le 
e(xx)^*e=ex^*x^*e$, and $S$ is $*$-intra-regular. $\hfill\Box$ 
\smallskip

\noindent{\bf Proposition 27.} {\it Let S be an involution 
$*$-intra-regular poe-semigroup. Then, for any $x\in S$, we 
have$$ex^*e\in (x^*)_{\cal N} \mbox { and } ex^*e\ge y \mbox { for 
every } y\in (x)_{\cal N}.$$}{\bf Proof.} Let $x\in S$. Since $e\ge 
x^*\in N(x^*)$, we have $e\in N(x^*)$. Since $e, x^*\in N(x^*)$, we 
have $ex^*e\in N(x^*)$, then $N(ex^*e)\subseteq N(x^*)$. Since 
$ex^*e\in N(ex^*e)$, we have $x^*\in N(ex^*e)$, then $N(x^*)\subseteq 
N(ex^*e)$. Thus we have $N(ex^*e)=N(x^*)$, so $ex^*e\in (x^*)_{\cal 
N}$. Let now $y\in (x)_{\cal N}$. Since $x\in N(x)=N(y)$, by Theorem 
26, we have $y\le ex^*e$. $\hfill\Box${\small

\end{document}